\documentclass[12pt]{article}

\usepackage{amssymb}
\usepackage{amsmath}
\usepackage{cite}
\usepackage[arrow, matrix, curve]{xy}

\begin{document}

\renewcommand{\citeleft}{{\rm [}}
\renewcommand{\citeright}{{\rm ]}}
\renewcommand{\citepunct}{{\rm,\ }}
\renewcommand{\citemid}{{\rm,\ }}

\newcounter{abschnitt}
\newtheorem{satz}{Theorem}
\newtheorem{coro}[satz]{Corollary}
\newtheorem{theorem}{Theorem}[abschnitt]
\newtheorem{koro}[theorem]{Corollary}
\newtheorem{prop}[theorem]{Proposition}
\newtheorem{lem}[theorem]{Lemma}
\newtheorem{expls}[theorem]{Examples}

\newcommand{\mres}{\mathbin{\vrule height 1.6ex depth 0pt width 0.11ex \vrule height 0.11ex depth 0pt width 1ex}}

\renewenvironment{quote}{\list{}{\leftmargin=0.62in\rightmargin=0.62in}\item[]}{\endlist}

\newcounter{saveeqn}
\newcommand{\alpheqn}{\setcounter{saveeqn}{\value{abschnitt}}
\renewcommand{\theequation}{\mbox{\arabic{saveeqn}.\arabic{equation}}}}
\newcommand{\reseteqn}{\setcounter{equation}{0}
\renewcommand{\theequation}{\arabic{equation}}}

\hyphenpenalty=9000

\sloppy

\phantom{a}

\vspace{-1.4cm}

\begin{center}
\begin{Large} {\bf Af\kern0pt f\kern0pt ine vs.\ Euclidean Isoperimetric Inequalities} \\[0.7cm] \end{Large}

\begin{large} Christoph Haberl and Franz E. Schuster\end{large}
\end{center}

\vspace{-0.8cm}

\begin{quote}
\footnotesize{ \vskip 1cm \noindent {\bf Abstract.}
It is shown that every even, zonal measure on the Euclidean unit sphere gives rise
to an isoperimetric inequality for sets of finite perimeter which directly implies the classical Euclidean
isoperimetric inequality. The strongest member of this large family of inequalities is shown to be
the only affine invariant one among them -- the Petty projection inequality. As an application, a
family of sharp Sobolev inequalities for functions of bounded variation is obtained, each of which is stronger than the classical Sobolev inequality.
Moreover, corresponding families of $L_p$ isoperimetric and Sobolev type inequalities are also established.}
\end{quote}

\vspace{0.6cm}

\centerline{\large{\bf{ \setcounter{abschnitt}{1}
\arabic{abschnitt}. Introduction}}}

\alpheqn

\vspace{0.6cm}

Over the last two decades several important affine isoperimetric inequalities, comparing geometric functionals which are
invariant under volume preserving affine or linear transformations, have been established (see, e.g., \textbf{\cite{boeroezky2013, Campi:Gronchi02a, habschu09, LYZ2000a, LYZ2010a}} or
the books \textbf{\cite{gardner2ed, schneider93}} for more information).
Despite that, it is still a common misbelief among geometers that Euclidean inequalities, that is, inequalities
for functionals invariant merely under rigid motions, are stronger than their affine counterparts.
One reason for this misconception might be that, so far, there are only a few explicit examples of Euclidean isoperimetric inequalities
that have been significantly improved by an affine invariant one. The best known and most important instance is
\emph{the} classical Euclidean isoperimetric inequality which is considerably strengthened by Petty's projection inequality \textbf{\cite{petty67}}.
Proved in the early 1970s, the latter can be seen as an integral geometric counterpart to the classical affine isoperimetric inequality from affine differential geometry.
In its original form, Petty's inequality states that among convex bodies of given volume, ellipsoids are precisely those whose polar projection bodies (see Section 2 for definitions) have maximal volume.
Subsequently, the Petty projection inequality has been generalized first to compact domains with smooth boundary by Zhang \textbf{\cite{zhang99}} and, more recently,
to sets of finite perimeter by Wang \textbf{\cite{wang12}}.

Even though Petty published his projection inequality in conference proceedings only, its tremendous impact to convex geometric analysis could hardly be overstated.
For example, it is the geometric core of the affine Sobolev-Zhang inequality \textbf{\cite{zhang99}} which strengthens the classical sharp Euclidean Sobolev inequality.
The $L_p$ version of Petty's inequality by Lutwak, Yang, and Zhang \textbf{\cite{LYZ2000a}} and its Orlicz extension by the same authors \textbf{\cite{LYZ2010a}} both marked landmark results
in the evolution of the Brunn-Minkowski theory first towards an $L_p$ theory and, more recently, towards an Orlicz theory of convex bodies.

\pagebreak

The special role of the Petty projection inequality has been further illuminated when Ludwig \textbf{\cite{ludwig02}} demonstrated that the projection body operator is the unique
continuous and affinely contravariant map on convex bodies which is a Minkowski valuation, that is, a finitely additive map with respect to Minkowski addition.
The notion of scalar valued valuations has long been an integral part of convex and discrete geometry (see, e.g., \textbf{\cite[\textnormal{Chapter 6}]{schneider93}} or \textbf{\cite{Alesker99, Alesker01, bernigfu10, habparap14, centro, parapwann}} for more recent results). The line of research concerned with characterizing Minkowski valuations is of newer vintage. Following Ludwig's seminal work,
a more or less complete picture on Minkowski valuations compatible with non-degenerate linear transformations was developed,
showing that they often form convex cones generated by finitely many maps (see \textbf{\cite{haberl11, Ludwig:Minkowski, Ludwig10a, SchuWan11, wannerer10}}). This is in stark contrast to the case of Minkowski valuations intertwining only rigid motions which form an infinite dimensional cone. Nonetheless, also here substantial inroads towards a complete classification have been made (see \textbf{\cite{kiderlen05, schneider74, Schu09, SchuWan13, SchuWan16}}), making it possible to extend affine inequalities for the projection body operator by Lutwak \textbf{\cite{lutwak93}} to a much larger class of Minkowski valuations (see \textbf{\cite{abardber11, ABS2011, BPSW2014, parapschu, Schu06a, Schu09}}). This raised the natural problem whether the Petty projection inequality also holds in greater generality.

\vspace{0.2cm}

In this article, we associate to every even, zonal measure on the Euclidean unit \linebreak sphere a unique Minkowski valuation intertwining rigid motions and prove that each \linebreak of these
operators gives rise to a sharp isoperimetric inequality for sets of finite perimeter.
This large family of inequalities has not only both the classical Euclidean isoperimetric inequality and the Petty projection inequality as special cases, but also
an inequality for the volume of polar mean section operators conjectured by Maresch and the second author \textbf{\cite{mareschschu}}.
Moreover, we show that each of these new inequalities strengthens and directly implies the Euclidean
isoperimetric inequality and we identify the strongest one. This turns out to be the only affine invariant inequality among them, the Petty projection inequality,
underlining how powerful affine inequalities truly are compared to their counterparts from Euclidean geometry.

Next, we show that each of our isoperimetric inequalities for sets of finite perimeter is equivalent to a Sobolev type inequality for functions of bounded variation.
This yields a family of analytic inequalities interpolating between the classical Sobolev inequality of Federer-Fleming \textbf{\cite{fedflem}} and Maz'ya \textbf{\cite{mazya}} and
the affine Sobolev-Zhang inequality \textbf{\cite{zhang99}}. As in the geometric setting, the only affine invariant inequality among them turns out to be the strongest one.

Finally, we obtain an extension of the $L_p$ Petty projection inequality of Lutwak, Yang, and Zhang \textbf{\cite{LYZ2000a}}
to a family of $L_p$ Minkowski valuations parameterized by even, zonal measures on the unit sphere. These $L_p$ isoperimetric inequalities are then used
to establish a family of new $L_p$ Sobolev type inequalities that interpolate between the sharp $L_p$ Sobolev inequality of Aubin \textbf{\cite{aubin76}} and Talenti \textbf{\cite{talenti76}} and the
affine $L_p$ Sobolev inequality of Lutwak, Yang, and Zhang \textbf{\cite{LYZ2002}}. As before, the only affine invariant member of the respective family of inequalities is the strongest
one.

\pagebreak

\vspace{1cm}

\centerline{\large{\bf{ \setcounter{abschnitt}{2}
\arabic{abschnitt}. Statement of principal results}}}

\reseteqn \alpheqn \setcounter{theorem}{0}

\vspace{0.6cm}

The setting for this article is Euclidean space $\mathbb{R}^n$, where we assume throughout that $n \geq 3$.
In addition to its denoting absolute value and the standard Euclidean norm on $\mathbb{R}^n$, we also write $|\cdot|$ for $k$-dimensional Hausdorff measure $\mathcal{H}^k$ in $\mathbb{R}^n$ for the
appropriate $k \in \{1, \ldots, n\}$. For example, if $K \subseteq \mathbb{R}^n$ is a convex body (that is, a compact, convex set) with nonempty interior, then $|K|$ denotes its volume and $|\partial K|$ its surface area. Similarly, $dx$ will always mean $d\mathcal{H}^k(x)$ for the appropriate $k$.

A measure on the Euclidean unit sphere $\mathbb{S}^{n-1}$ of $\mathbb{R}^n$ is to be understood to mean a \emph{non-negative} and \emph{non-trivial} finite Borel measure
on $\mathbb{S}^{n-1}$. We use $\bar{e} \in \mathbb{S}^{n-1}$ to denote an arbitrary but fixed point (the pole) of
the sphere and we write $\mathrm{SO}(n-1)$ for the stabilizer of $\bar{e}$ in $\mathrm{SO}(n)$. A measure on $\mathbb{S}^{n-1}$ is
said to be \emph{even} if it assigns the same value to antipodal sets and it is called \emph{zonal} if it is $\mathrm{SO}(n-1)$ invariant.

Let $\mathcal{K}^n$ be the space of convex bodies in $\mathbb{R}^n$ endowed with the Hausdorff metric
and recall that each $K \in \mathcal{K}^n$ is uniquely determined by its support function $h(K,u) = \max\{u \cdot x: x \in K\}$ for $u \in \mathbb{S}^{n-1}$.
If $K$ contains the origin, then $K^{\circ} = \{x \in \mathbb{R}^n: x \cdot y \leq 1 \mbox{ for all } y \in K\}$ is the \emph{polar body} of $K$.

For a set of finite perimeter $L \subseteq \mathbb{R}^n$, we denote by $\partial^* L$ the \emph{reduced boundary} of $L$
and we write $\nu_L$ for the (measure theoretic) \emph{outer unit normal vector} field to $L$ (see Section 3 for detailed definitions).
The \emph{projection body} of $L$ is the convex body $\Pi L \in \mathcal{K}^n$ with support function
\begin{equation} \label{defpi}
h(\Pi L,u) = \frac{1}{2} \int_{\partial^*L} |u \cdot \nu_L(x) |\,dx, \qquad u \in \mathbb{S}^{n-1}.
\end{equation}
This extension of Minkowski's classical notion of the projection body of a convex body was first given by Wang \textbf{\cite{wang12}}, who generalized
a definition of Zhang \textbf{\cite{zhang99}} for compact sets with piecewise $C^1$ boundary.
Wang \textbf{\cite{wang12}} also established the \emph{Petty projection inequality} for sets of finite perimeter (generalizing the original result of Petty \textbf{\cite{petty67}} for convex bodies and
a version for compact sets with piecewise $C^1$ boundary by Zhang \textbf{\cite{zhang99}}): A set of finite perimeter $L \subseteq \mathbb{R}^n$ is a maximizer of
the volume product $|\Pi^{\circ}L||L|^{n-1}$ if and only if $L$ is an ellipsoid up to a Lebesgue null set. \linebreak
(Here and in the following, we write $\Pi^{\circ}L$ instead of $(\Pi L)^{\circ}$.)

For $p \geq 1$, each even measure $\mu$ on $\mathbb{S}^{n-1}$ generates an origin-symmetric convex body $Z_p^{\mu} \in \mathcal{K}^n$ (which is uniquely determined when $p$ is not an even integer) by
\begin{equation} \label{zonoidmeas}
h(Z^{\mu}_p,u)^p = \int_{\mathbb{S}^{n-1}} |u \cdot v|^p\,d\mu(v), \qquad u \in \mathbb{S}^{n-1}.
\end{equation}
The bodies obtained in this way constitute the class of origin-symmetric \emph{$L_p$~zonoids}, which also arise naturally in various other contexts (see, e.g., \textbf{\cite[\textnormal{Chapter~3.5}]{schneider93}}). When $p = 1$, $L_p$ zonoids are just called zonoids and we write $Z^\mu$ rather than $Z_1^{\mu}$. If $\mu$ is zonal, we use $Z^{\mu}_p(\bar{e})$ instead of $Z_p^{\mu}$ to indicate the bodies axis of symmetry. The rotated copy of $Z^{\mu}_p(\bar{e})$ whose axis of symmetry is $u \in \mathbb{S}^{n-1}$ is denoted by $Z_p^{\mu}(u)$.

We are now in a position to define the main objects of our investigations.

\vspace{0.3cm}

\noindent {\bf Definition} \emph{Suppose that $\mu$ is an even, zonal measure on $\mathbb{S}^{n-1}$ and let $Z^{\mu}(\bar{e})$ be the zonoid generated by $\mu$.
For a set of finite perimeter $L \subseteq \mathbb{R}^n$, we define the convex body $\Phi^{\mu}L \in \mathcal{K}^n$ by}
\begin{equation} \label{defphimu}
h(\Phi^{\mu}L,u) = \int_{\partial^* L} h(Z^{\mu}(u),\nu_L(x))\,dx, \qquad u \in \mathbb{S}^{n-1}.
\end{equation}

\vspace{0.3cm}

We shall see in the next section that (\ref{defphimu}) indeed defines a support function for every set of finite perimeter.
Before we discuss the origins and basic properties of the operators $\Phi^{\mu}$, let us look at some special cases of their construction.

\vspace{0.2cm}

\noindent {\bf Examples:}

\begin{enumerate}
\item[(a)] If $\mu$ is a multiple of spherical Lebesgue measure, then, by (\ref{zonoidmeas}), $Z^{\mu}$ is a dilate of the Euclidean unit ball $B^n$ and, thus, $h(Z^{\mu},\cdot)$ is
proportional to the standard Euclidean norm on $\mathbb{R}^n$. Consequently, by (\ref{defphimu}), we have in this case
\[\Phi^{\mu}L \cong |\partial^* L| B^n  \]
for every set of finite perimeter $L \subseteq \mathbb{R}^n$, where $|\partial^* L|$ is the perimeter of $L$ (see Section 3). In particular, if $L$ is a convex body or has piecewise smooth boundary, then $|\partial^* L|$ is just the usual surface area of $L$.

\item[(b)] If $\mu$ is discrete, then, since $\mu$ is even and zonal, it must be a multiple of the sum of two Dirac measures $\delta_{\bar{e}} + \delta_{-\bar{e}}$.
Thus, by (\ref{zonoidmeas}), $Z^{\mu}(\bar{e})$ is a dilate of the segment $[-\bar{e},\bar{e}]$ and $h(Z^{\mu}(u),v) \cong |u \cdot v|$. Hence, by (\ref{defpi}) and (\ref{defphimu}), we obtain
\[\Phi^{\mu}L \cong \Pi L  \]
for every set of finite perimeter $L \subseteq \mathbb{R}^n$.

\item[(c)] If $\mu$ is concentrated on the equator of $\mathbb{S}^{n-1}$ perpendicular to the pole $\bar{e}$,
then (\ref{zonoidmeas}) implies that $h(Z^{\mu}(u),v) \cong \sqrt{1-(u\cdot v)^2}$ and $Z^{\mu}(\bar{e})$ is a dilate of the disc $B^n \cap \bar{e}^{\bot}$.
Thus, by (\ref{defphimu}), we have for every set of finite perimeter $L \subseteq \mathbb{R}^n$,
\[\Phi^{\mu}L \cong \mathrm{M}_2^+ L,  \]
where $\mathrm{M}_2^+$ denotes the even part of the \emph{second mean section operator}. Mean section operators were introduced in 1992 by Goodey and Weil \textbf{\cite{GooWei92}} and have recently
become the focus of increased interest (see, e.g., \textbf{\cite{GooWei12, goodeyweil4, mareschschu, SchuWan13,SchuWan16}}).

\end{enumerate}

The operators $\Phi^{\mu}$ first appeared about a decade ago in the theory of Minkowski valuations (see \textbf{\cite{Schu06a}}).
In general, a \emph{valuation} on the space $\mathcal{K}^n$ is a map $\Psi: \mathcal{K}^n \rightarrow \mathcal{A}$ with values in an Abelian semigroup $\mathcal{A}$ such that
\[\Psi(K) + \Psi(L) = \Psi(K \cup L) + \Psi(K \cap L)   \]
whenever $K \cup L$ is convex. When $\mathcal{A} = \mathcal{K}^n$ and addition on $\mathcal{K}^n$ is the usual Minkowski addition, then $\Psi$ is called a Minkowski valuation.

\pagebreak

It is not difficult to see that the restriction of $\Phi^{\mu}$ to convex bodies defines a continuous Minkowski valuation $\Phi^{\mu}: \mathcal{K}^n \rightarrow \mathcal{K}^n$
which is translation invariant and commutes with $\mathrm{SO}(n)$. Moreover, it follows from a characterization result of Ludwig \textbf{\cite{Ludwig:Minkowski}} that $\Phi^{\mu}$
is compatible with $\mathrm{SL}(n)$ (more precisely, $\mathrm{SL}(n)$ contravariant) if and only if $\mu$ is discrete, that is, if and only if $\Phi^{\mu}$ is a multiple
of $\Pi$. It was also first asked by Ludwig whether an isoperimetric inequality similar to Petty's projection inequality holds for the operators $\Phi^{\mu}$.
With our first main result, we answer this question in the affirmative. (Throughout we use the convention that $\infty \cdot 0 = 0$.)

\begin{satz} \label{mainthm1} Suppose that $\mu$ is an even, zonal measure on $\mathbb{S}^{n-1}$. Among sets of finite perimeter $L \subseteq \mathbb{R}^n$ the volume product
\[|\Phi^{\mu,\circ}L||L|^{n-1}   \]
is maximized by Euclidean balls. If $\mu$ is not discrete, then Euclidean balls are the only maximizers up to Lebesgue null sets. If $\mu$ is discrete, then
$L$ is a maximizer if and only if it is an ellipsoid up to a Lebesgue null set.
\end{satz}

Note that by the examples above the Euclidean isoperimetric inequality and the Petty projection inequality are both special cases of Theorem \ref{mainthm1}. Moreover, if
$\mu$ is concentrated on the equator of $\mathbb{S}^{n-1}$ perpendicular to $\bar{e}$, then Theorem \ref{mainthm1} confirms a conjecture of Maresch and the second author \textbf{\cite{mareschschu}}
about the maximizers of the volume product $|\mathrm{M}_2^{+,\circ} L||L|^{n-1}$.

It was first pointed out by Lutwak \textbf{\cite{lutwak84}} that the Petty projection inequality is significantly stronger than the classical isoperimetric
inequality. As our next result shows, each new inequality of Theorem \ref{mainthm1} is, in fact, stronger and directly implies the isoperimetric inequality. Moreover,
the Petty projection inequality -- the only affine invariant one among them -- is the strongest one.

\begin{satz} \label{mainthm2} If $\mu$ is an even, zonal measure on $\mathbb{S}^{n-1}$ such that $\mu(\mathbb{S}^{n-1}) = \frac{1}{2}$ and $L \subseteq \mathbb{R}^n$ is a set of finite perimeter with nonempty interior, then
\[\frac{n^n\omega_n^{n+1}}{\omega_{n-1}^n}|\partial^*L|^{-n} \leq |\Phi^{\mu,\circ}L| \leq |\Pi^\circ L|.  \]
There is equality in the left hand inequality if and only if $\Phi^{\mu}L$ is a Euclidean ball.
There is equality in the right hand inequality if and only if $\mu$ is discrete or $\Pi L$ is a Euclidean ball.
\end{satz}

Here, $\omega_k=\pi^{k/2}/\Gamma(1+k/2)$ is the $k$-dimensional Lebesgue measure of the Euclidean unit ball in $\mathbb{R}^k$ and the normalization $\mu(\mathbb{S}^{n-1}) = \frac{1}{2}$ is chosen such that
$\Phi^{\mu}B^n = \Pi B^n = \omega_{n-1} B^n$ for all $\mu$. Note that the right inequality of Theorem \ref{mainthm2} combined with the Petty projection inequality implies Theorem \ref{mainthm1} (up to a short argument to obtain
the characterization of the extremizers). In Section 5, we will therefore first prove Theorem \ref{mainthm2} and deduce Theorem \ref{mainthm1} as a consequence.

\pagebreak

It is a classical fact (see, e.g., \textbf{\cite{evansgariepy}}) that the functional form of the isoperimetric inequality is the \emph{Sobolev inequality} on the space $BV(\mathbb{R}^n)$ of functions of bounded variation. It states that for every $f \in BV(\mathbb{R}^n)$,
\begin{equation} \label{L1Sob}
\|Df\| \geq n\omega_n^{1/n} \|f\|_{\frac{n}{n-1}},
\end{equation}
where the vector valued Radon measure $Df$ is the weak gradient of $f$, $\|Df\|$ denotes its total variation in $\mathbb{R}^n$ (see Section 3), and $\|f\|_p$ denotes the usual $L_p$ norm of $f$ in $\mathbb{R}^n$.
Equality holds in (\ref{L1Sob}) if and only if $f$ is a multiple of the characteristic function of a Euclidean ball in $\mathbb{R}^n$.

An analytic version of the Petty projection inequality was first established by Zhang \textbf{\cite{zhang99}} for compactly supported $C^1$ functions. It became known as the \emph{affine
Sobolev-Zhang inequality} and was, more recently, extended to $BV(\mathbb{R}^n)$ by Wang~\textbf{\cite{wang12}}. This extended inequality states that
for every $f \in BV(\mathbb{R}^n)$,
\begin{equation} \label{SobZhang}
\left ( \int_{\mathbb{S}^{n-1}} \left ( \int_{\mathbb{R}^n} |u \cdot \sigma_f| \,d|Df| \right )^{-n} du \right )^{-1/n} \geq \frac{2\omega_{n-1}}{n^{1/n}\omega_n}\, \|f\|_{\frac{n}{n-1}}
\end{equation}
with equality if and only if $f$ is a multiple of the characteristic function of an ellipsoid in $\mathbb{R}^n$.
Here, we write $|Df|$ for the variation measure of $Df$ and $\sigma_f$ is the Radon-Nikodym derivative of $Df$ with respect to $|Df|$.
Using H\"older's inequality, it follows that the affine invariant inequality (\ref{SobZhang}) directly implies the classical inequality (\ref{L1Sob}) which is merely invariant under rigid motions.

The proof of (\ref{SobZhang}) given by Wang \textbf{\cite{wang12}} was based on an approach towards the affine
Sobolev-Zhang inequality for compactly supported $C^1$ functions developed by Lutwak, Yang, and Zhang \textbf{\cite{LYZ2006}}.
Using the same ideas, we obtain the following family of sharp Sobolev type inequalities which can be seen as the functional form of the isoperimetric inequalities from Theorem \ref{mainthm1}.

\begin{satz} \label{mainthm3} Suppose that $\mu$ is an even, zonal measure on $\mathbb{S}^{n-1}$ and let $Z^{\mu}(\bar{e})$ be the zonoid generated by $\mu$.
Then, for every $f \in BV(\mathbb{R}^n)$,
\[\left ( \int_{\mathbb{S}^{n-1}} \left ( \int_{\mathbb{R}^n} \|\sigma_f\|_{Z^{\mu}(u)^\circ}\,d|Df| \right )^{-n} du \right )^{-1/n} \geq \frac{2\omega_{n-1}\mu(\mathbb{S}^{n-1})}{n^{1/n}\omega_n}\, \|f\|_{\frac{n}{n-1}}.\]
If $\mu$ is not discrete, there is equality if and only if $f$ is a multiple of the characteristic function of a Euclidean ball in $\mathbb{R}^n$.
If $\mu$ is discrete, then equality holds if and only if $f$ is a multiple of the characteristic function of an ellipsoid in $\mathbb{R}^n$.
\end{satz}

Here, $\|\cdot\|_{Z^{\mu}(u)^\circ} = h(Z^{\mu}(u),\cdot)$ denotes the seminorm on $\mathbb{R}^n$ (which is actually a norm when $\mu$ is not discrete) with unit ball $Z^{\mu}(u)^\circ$.
We will see in Section 6 how Theorem~\ref{mainthm2} implies that each of the new analytic inequalities from Theorem~\ref{mainthm3} is stronger than the Sobolev inequality (\ref{L1Sob})
and, moreover, that the only affine invariant inequality from Theorem \ref{mainthm3} -- the affine Sobolev-Zhang inequality -- is the strongest one.

\pagebreak

A stronger affine version of the sharp $L_p$ Sobolev inequality of Aubin \textbf{\cite{aubin76}} and Talenti \textbf{\cite{talenti76}} was
established by Lutwak, Yang, and Zhang \textbf{\cite{LYZ2002}}. If $W^{1,p}(\mathbb{R}^n)$ denotes the Sobolev space of real-valued $L_p$ functions on $\mathbb{R}^n$
with weak $L_p$ partial derivatives, then the \emph{affine $L_p$ Sobolev inequality} states that for $1 < p < n$ and $f \in W^{1,p}(\mathbb{R}^n)$,
\begin{equation}\label{affLpSob}
\left ( \int_{\mathbb{S}^{n-1}} \left ( \int_{\mathbb{R}^n} |\nabla f(x) \cdot u|^p\, dx \right )^{-n/p} du \right )^{-1/n} \geq c_{n,p}\, \|f\|_{p^*},
\end{equation}
where $p^*=np/(n-p)$. The optimal constants $c_{n,p}$ in (\ref{affLpSob}) were explicitly computed in \textbf{\cite{LYZ2002}}.
It was later shown by Wang \textbf{\cite{wang13}} (see also \textbf{\cite{nguyen}} for a different approach) that under mild additional technical assumptions on $f$, equality holds in (\ref{affLpSob})
if and only if $f$ coincides up to translation a.e.\ on $\mathbb{R}^n$ with its convex symmetrization $f^E$ with respect to an origin-symmetric ellipsoid $E$ in $\mathbb{R}^n$ (cf.\ \textbf{\cite{wang13}} or \textbf{\cite{nguyen}} for definitions).

We establish an $L_p$ version of Theorem \ref{mainthm3} which shows that also the classical and the affine $L_p$ Sobolev inequalities are members of a larger family of analytic
inequalities parametrized by even, zonal measures on $\mathbb{S}^{n-1}$.

\begin{satz} \label{mainthm4}  Suppose that $1 < p < n$ and that $\mu$ is an even, zonal measure on $\mathbb{S}^{n-1}$. If $Z^{\mu}_p(\bar{e})$ is the $L_p$ zonoid generated by $\mu$,
then, for every $f \in W^{1,p}(\mathbb{R}^n)$,
\begin{equation} \label{HSgenlpsob}
\left ( \int_{\mathbb{S}^{n-1}} \left ( \int_{\mathbb{R}^n} \|\nabla f(x)\|^p_{Z^{\mu}_p(u)^\circ}\,dx \right )^{-n/p} du \right )^{-1/n} \geq c_{n,p}\mu(\mathbb{S}^{n-1})^{1/p} \|f\|_{p^*}.
\end{equation}
\end{satz}

\vspace{0.2cm}

If $\mu$ is discrete, then inequality (\ref{HSgenlpsob}) reduces to (\ref{affLpSob}) and if $\mu$ is proportional to spherical Lebesgue measure, then (\ref{HSgenlpsob}) becomes the classical $L_p$ Sobolev inequality of Aubin and Talenti. Using either the approach from \textbf{\cite{wang13}} or that of \textbf{\cite{nguyen}}, it is possible to show that if $\mu$ is not discrete, then
equality holds in (\ref{HSgenlpsob}) if and only if $f$ coincides up to translation a.e.\ on $\mathbb{R}^n$ with its symmetric rearrangement $f^{\mbox{\begin{tiny}$\bigstar$\end{tiny}}}$. However, since this requires additional tools and techniques different from the ones applied to prove Theorem \ref{mainthm3}, it will be explored in a later publication (see \textbf{\cite{KniefaczSchu}}).

The geometric inequalities needed to establish Theorem \ref{mainthm4} are the content of an $L_p$ version of Theorem \ref{mainthm1} which generalizes the $L_p$ Petty projection inequality of Lutwak, Yang, and Zhang \textbf{\cite{LYZ2000a}} (see \textbf{\cite{Campi:Gronchi02a}} and \textbf{\cite{boeroezky2013}} for alternate approaches). It is important to note that while the geometry behind the classical $L_p$ Sobolev inequality (that is, the isoperimetric inequality) is the same for all $p$, the $L_p$ isoperimetric inequalities behind Theorem \ref{mainthm4} are different for different $p$.
We establish this family of inequalities in Section 5 and show that the strongest member is again the only affine invariant one among them -- the $L_p$ Petty projection inequality. This also implies that the affine $L_p$ Sobolev inequality is the strongest inequality among the family of Sobolev inequalities (\ref{HSgenlpsob}).

\vspace{1cm}

\centerline{\large{\bf{ \setcounter{abschnitt}{3}
\arabic{abschnitt}. Background material}}}

\reseteqn \alpheqn \setcounter{theorem}{0}

\vspace{0.6cm}

In this section we recall basic definitions and facts about functions of bounded variation and sets of finite perimeter as well as some notions and results from the
$L_p$ Brunn-Minkowski theory of convex bodies required for the proofs of our main results. As general references for this material we refer to the classic book by Evans and Gariepy \textbf{\cite{evansgariepy}} and the recent monograph \textbf{\cite{schneider93}} by Schneider.

\vspace{0.2cm}

Let $C_c^k(\mathbb{R}^n)$ denote the space of all compactly supported $C^k$ functions on $\mathbb{R}^n$.
A function $f \in L^1(\mathbb{R}^n)$ is called a \emph{function of bounded variation} on $\mathbb{R}^n$ if for every $1 \leq i \leq n$, there exists a finite signed Radon measure $D_if$ on $\mathbb{R}^n$ such that
\begin{equation} \label{bvradonmeas}
\int_{\mathbb{R}^n} f\,\frac{\partial g}{\partial x_i}\,dx = - \int_{\mathbb{R}^n} g\,dD_if
\end{equation}
for all $g \in C_c^1(\mathbb{R}^n)$. We denote the space of all functions of bounded variation on $\mathbb{R}^n$ by $BV(\mathbb{R}^n)$. (Note that, as usual, two functions in $BV(\mathbb{R}^n)$ that coincide almost everywhere with respect to Lebesgue measure are considered to be the same.)

For the vector valued Radon measure $Df = (D_1f,\ldots,D_nf)$ on $\mathbb{R}^n$, we define its (Euclidean) \emph{variation} as the non-negative Radon measure $|Df|$ whose value at a Borel set $E \subseteq \mathbb{R}^n$ is given by
\[|Df|(E) = \sup_{\pi} \sum_{A \in \pi} |Df(A)|,  \]
where the supremum is taken over all partitions $\pi$ of $E$ into a countable number of disjoint measurable subsets. For $f \in BV(\mathbb{R}^n)$, we denote the Radon-Nikodym derivative of $Df$ with respect to $|Df|$ by $\sigma_f$. Hence, it follows from (\ref{bvradonmeas}), that the vector field $\sigma_f$ satisfies
\begin{equation} \label{intsatzbv1}
\int_{\mathbb{R}^n} f\,\mathrm{div}\phi\,dx = - \int_{\mathbb{R}^n} \phi\cdot dDf = - \int_{\mathbb{R}^n} \phi\cdot \sigma_f\,d|Df|
\end{equation}
for all $\phi \in C_c^1(\mathbb{R}^n,\mathbb{R}^n)$, the space of all continuously differentiable vector fields on $\mathbb{R}^n$ with compact support.

The \emph{perimeter} of a (Lebesgue) measurable set $L \subseteq \mathbb{R}^n$ is defined by
\[P(L)=  \sup \left \{ \int_L \mathrm{div} \phi\,dx: \phi \in C_c^1(\mathbb{R}^n,\mathbb{R}^n), \|\phi\|_\infty \leq 1   \right \}.\]
If $P(L)<\infty$, then $L$ is called a \emph{set of finite perimeter}. Using Riesz's representation theorem and (\ref{intsatzbv1}), it is not difficult to show that $L \subseteq \mathbb{R}^n$ has finite perimeter if and only if $1_L \in BV(\mathbb{R}^n)$.

For $x \in \mathbb{R}^n$, we denote by $B_r(x)$ the Euclidean ball with center $x$ and radius $r$. \linebreak
A (topological) boundary point $x$ of a set of finite perimeter $L \subseteq \mathbb{R}^n$ is said to belong to the \emph{reduced boundary} $\partial^*L$ of $L$ if
\[\lim_{r \rightarrow 0^+} \frac{D1_L(B_r(x))}{|D1_L|(B_r(x))}  \]
exists and belongs to $\mathbb{S}^{n-1}$. In this case, we call the limit $\nu_L(x) = -\sigma_{1_L}(x)$ the \emph{measure-theoretic outer unit normal} to $L$ at $x$. It follows from De Giorgi's structure theorem that
\begin{equation} \label{Dvonchar}
D1_L = \nu_Ldx \mres \partial^*L \qquad \mbox{and} \qquad |D1_L|(\mathbb{R}^n) = |\partial^*L|.
\end{equation}

Note that if $L$ is a compact set with $C^1$ boundary, then $\partial^*L = \partial L$ and the measure-theoretic outer unit normals coincide with the usual outer unit normals.

\vspace{0.2cm}

Recall that $\mathcal{K}^n$ denotes the space of convex bodies in $\mathbb{R}^n$ and that each $K \in \mathcal{K}^n$ is uniquely determined by its
support function $h(K,u) = \max\{u \cdot x: x \in K\}$, $u \in \mathbb{S}^{n-1}$. It follows from this definition that for each $\vartheta \in \mathrm{SO}(n)$, we have
\begin{equation} \label{sptfctcomprot}
h(\vartheta K,u) = h(K,\vartheta^{-1}u), \qquad u \in \mathbb{S}^{n-1}.
\end{equation}
The \emph{mean width} of a convex body $K \in \mathcal{K}^n$ is defined by
\[w(K) = \frac{2}{n\omega_n} \int_{\mathbb{S}^{n-1}} h(K,u)\,du.  \]
Using (\ref{zonoidmeas}), Fubini's theorem, and the Cauchy projection formula (see \ref{pisurfmeas} below), we can calculate the mean width of a zonoid $Z^{\mu} \in \mathcal{K}^n$, generated by the even measure $\mu$ on $\mathbb{S}^{n-1}$, as follows:
\begin{equation} \label{meanwidthzonoid}
w(Z^{\mu}) = \frac{2}{n\omega_n} \int_{\mathbb{S}^{n-1}} \int_{\mathbb{S}^{n-1}} |u \cdot v|\,du\,d\mu(v) = \frac{4\omega_{n-1}}{n\omega_n}\mu(\mathbb{S}^{n-1}).
\end{equation}
Urysohn's inequality states that for every $K \in \mathcal{K}^n$,
\begin{equation} \label{urysohn}
\frac{w(K)}{2} \geq \frac{|K|^{1/n}}{\omega_n^{1/n}}
\end{equation}
with equality if and only if $K$ is a ball.

A convex body $K \in \mathcal{K}^n$ containing the origin is also uniquely determined by its radial function $\rho(K,u) = \max\{\lambda \geq 0: \lambda u \in K\}$, $u \in \mathbb{S}^{n-1}$.
Note that if the origin is an interior point of $K$, then
\[\rho(K^{\circ},\cdot) = \frac{1}{h(K,\cdot)}, \qquad h(K^{\circ},\cdot) = \frac{1}{\rho(K,\cdot)}, \qquad h(K,\cdot) = \|\cdot\|_{K^{\circ}},  \]
where $K^{\circ} = \{x \in \mathbb{R}^n: x \cdot y \leq 1 \mbox{ for all } y \in K\}$ is the polar body of $K$. In particular, the polar coordinate formula for the volume of $K$ takes the form
\begin{equation} \label{polarcoordvol}
|K| = \frac{1}{n} \int_{\mathbb{S}^{n-1}} \rho(K,u)^n\,du = \frac{1}{n} \int_{\mathbb{S}^{n-1}} h(K^\circ,u)^{-n}\,du.
\end{equation}

\pagebreak

For $K, L \in \mathcal{K}^n$ and $t > 0$, the \emph{Minkowski combination} $K + tL \in \mathcal{K}^n$ can be defined by
\[h(K + tL,\cdot) = h(K,\cdot) + t\,h(L,\cdot).  \]
The \emph{surface area measure} $S(K,\cdot)$ of a convex body $K \in \mathcal{K}^n$ is the unique Borel measure on $\mathbb{S}^{n-1}$ such that
\[\lim_{t\rightarrow 0^+} \frac{|K + tL| - |K|}{t} = \int_{\mathbb{S}^{n-1}} h(L,u)\,dS(K,u)  \]
for each $L \in \mathcal{K}^n$. Alternatively, the surface area measure of $K$ can be defined by setting
\begin{equation} \label{surfareagauss}
\int_{\mathbb{S}^{n-1}} f(u)\,dS(K,u) = \int_{\partial' K} f(\nu_K(x))\,dx
\end{equation}
for each $f \in C(\mathbb{S}^{n-1})$. Here, the Gauss map $\nu_K: \partial'K \rightarrow \mathbb{S}^{n-1}$
is defined on the subset $\partial'K$ of those points of $\partial K$ that have a unique outer unit normal and is hence defined $\mathcal{H}^{n-1}$ a.e.\ on $\partial K$.

Note that for $\vartheta \in \mathrm{SO}(n)$, the measure $S(\vartheta K,\cdot)$ is just the pushforward of $S(K,\cdot)$ under the rotation $\vartheta$.
Another basic feature of the surface area measure is the valuation property, that is,
\begin{equation} \label{valpropsurfmeas}
S(K \cup L,\cdot) + S(K \cap L,\cdot) = S(K,\cdot) + S(L,\cdot)
\end{equation}
whenever $K \cup L \in \mathcal{K}^n$.

Minkowski defined the \emph{projection body} of a convex body $K \in \mathcal{K}^n$ as the unique convex body whose support function is given by
\begin{equation} \label{pisurfmeas}
h(\Pi K,u) = \mathrm{vol}_{n-1}(K|u^{\bot}) = \frac{1}{2} \int_{\mathbb{S}^{n-1}} |u \cdot v|\,dS(K,v), \qquad u \in \mathbb{S}^{n-1}.
\end{equation}
Here, the second equality is known as Cauchy's projection formula. Note that, by (\ref{valpropsurfmeas}), the projection body map is a Minkowski valuation.
Moreover, the following characterization result was obtained by Ludwig.

\begin{theorem} \label{ludchar} \emph{(Ludwig \textbf{\cite{Ludwig:Minkowski}})}
A map $\Psi: \mathcal{K}^n \rightarrow \mathcal{K}^n$ is a continuous and translation invariant Minkowski valuation such that $\Psi(AK) = A^{-\mathrm{T}}\Psi K$ for every $K \in \mathcal{K}^n$ and $A \in \mathrm{SL}(n)$ if and only if $\Psi = c\,\Pi$ for some $c \geq 0$.
\end{theorem}

It follows from (\ref{surfareagauss}) that not only does (\ref{defpi}) extend Minkowski's definition of the projection body of a convex body to sets of finite perimeter but also that
for every convex body $K \in \mathcal{K}^n$, definition (\ref{defphimu}) is equivalent to
\begin{equation} \label{phisurfmeas}
h(\Phi^{\mu} K,u) = \int_{\mathbb{S}^{n-1}} h(Z^{\mu}(u),v)\,dS(K,v), \qquad u \in \mathbb{S}^{n-1}.
\end{equation}
Since $Z^{\mu}(u)$ is a body of revolution, $h(Z^{\mu}(u),v)$ depends only on the value of $u \cdot v$ and, hence, $h(Z^{\mu}(u),v)=h(Z^{\mu}(v),u)$. This shows that (\ref{phisurfmeas}) and
(\ref{defphimu}) indeed define support functions of convex bodies.

\pagebreak

By (\ref{valpropsurfmeas}) and (\ref{phisurfmeas}), each map $\Phi^{\mu}$ is a Minkowski valuation. Moreover, using (\ref{sptfctcomprot}) and the fact that $\vartheta Z^{\mu}(u) = Z^{\mu}(\vartheta u)$, we see
that $\Phi^{\mu}(\vartheta K) = \vartheta \Phi^{\mu}K$ for every $\vartheta \in \mathrm{SO}(n)$ and $K \in \mathcal{K}^n$, that is, $\Phi^{\mu}$ is $\mathrm{SO}(n)$ equivariant.
In particular, since $(\vartheta K)^{\circ} = \vartheta K^{\circ}$, the quantities $|\Phi^{\mu,\circ}K|$ are rigid motion invariant. However, by Theorem \ref{ludchar},
$|\Phi^{\mu,\circ}K|$ is invariant under volume preserving affine transformations if and only if $\mu$ is discrete.

\vspace{0.2cm}

By using $L_p$ Minkowski combinations first introduced in the 1960s by Firey, Lutwak \textbf{\cite{lutwak93a, lutwak96}} showed that
the classical Brunn-Minkowski theory of convex bodies has a natural extension to a more general $L_p$ theory. Next we review some of the details of this
theory which are necessary for the proof of Theorem \ref{mainthm4}.

Suppose that $1 \leq p < \infty$ and that $K, L \in \mathcal{K}^n$ contain the origin in their interiors. For $t > 0$, the \emph{$L_p$ Minkowski combination} $K +_p t\cdot L \in \mathcal{K}^n$ can be defined by
\[h(K +_p t\cdot L,\cdot)^p = h(K,\cdot)^p + t\,h(L,\cdot)^p.  \]
The \emph{$L_p$ surface area measure} $S_p(K,\cdot)$ of $K$ is the unique Borel measure on $\mathbb{S}^{n-1}$ such that
\[\lim_{t\rightarrow 0^+} \frac{|K +_p t\cdot L| - |K|}{t} = \frac{1}{p} \int_{\mathbb{S}^{n-1}} h(L,u)^p\,dS_p(K,u)  \]
for each $L \in \mathcal{K}^n$ containing the origin in its interior. It was shown in \textbf{\cite{lutwak93a}} that $S_p(K,\cdot)$ is absolutely continuous with respect to
$S_1(K,\cdot) = S(K,\cdot)$ and that its Radon-Nikodym derivative is $h(K,\cdot)^{1-p}$.

The $L_p$ projection body of a convex body $K \in \mathcal{K}^n$ containing the origin in its interior was first defined in \textbf{\cite{LYZ2000a}} for $1 \leq p < \infty$ by
\begin{equation} \label{piplpsurfmeas}
h(\Pi_p K,u)^p = a_{n,p} \int_{\mathbb{S}^{n-1}} |u \cdot v|^p\,dS_p(K,v), \qquad u \in \mathbb{S}^{n-1},
\end{equation}
where the normalizing constant $a_{n,p}$ is chosen such that $\Pi_p B_1(0) = B_1(0)$ and, consequently, is given by
\[a_{n,p} = \frac{\Gamma\left ( \frac{n+p}{2} \right )}{2\pi^{(n-1)/2}\Gamma\left (\frac{p+1}{2}\right )}.  \]
Note that when $p = 1$, definition (\ref{piplpsurfmeas}) can be extended to all $K \in \mathcal{K}^n$ and that, in this case, we have
\begin{equation} \label{pi1pi}
\Pi_1 K = \omega_{n-1}^{-1}\Pi K.
\end{equation}

An $L_p$ analogue of Petty's projection inequality was first obtained in \textbf{\cite{LYZ2000a}} (see also \textbf{\cite{boeroezky2013, Campi:Gronchi02a, habschu09, Lin}} for different variants):
A convex body $K$ in $\mathbb{R}^n$ containing the origin in its interior is a maximizer of
the volume product $|\Pi_p^{\circ}K|^p|K|^{n-p}$ for $p > 1$ if and only if $K$ is an ellipsoid centered at the origin.

The geometric core in our proof of Theorem \ref{mainthm4} is an $L_p$ version of Theorem \ref{mainthm1}. To this end we introduce in the next section an $L_p$ extension of the operators $\Phi^{\mu}$ which is motivated by (\ref{piplpsurfmeas}). In Section 5, we then state and prove the required generalization of the $L_p$ Petty projection inequality.

\pagebreak

A few years after their first proof of the affine $L_p$ Sobolev inequality (\ref{affLpSob}), Lutwak, Yang, and Zhang \textbf{\cite{LYZ2006}} gave a new more conceptual proof of (\ref{affLpSob}) by associating to each $f \in W^{1,p}(\mathbb{R}^n)$ an origin-symmetric convex body $\langle f \rangle_p$. This convexification of a given Sobolev function is based on the solution of a functional analogue of the even $L_p$ Minkowski problem. The latter was first formulated in \textbf{\cite{lutwak93a}} and can be stated as follows: Given an even Borel measure $\nu$ on $\mathbb{S}^{n-1}$, does there exist an origin-symmetric convex body $K$ such that $\nu = S_p(K,\cdot)$? The case $p = 1$ is the classical (even) Minkowski problem which was solved by Minkowski, Aleksandrov, and Fenchel and Jessen (see, e.g., \textbf{\cite[\textnormal{Chapter 8.2}]{schneider93}}). Lutwak \textbf{\cite{lutwak93a}} gave an affirmative answer to the even $L_p$ Minkowski problem for all $p \neq n$; in \textbf{\cite{LYZ2004}} a \emph{volume-normalized} version of the $L_p$ Minkowski problem was introduced for which the case $p = n$ could be dealt with as well. Here, we state the solution to the functional volume-normalized even $L_p$ Minkowski problem obtained in \textbf{\cite{LYZ2006}} which is crucial in our proof of Theorem \ref{mainthm4}.

\begin{theorem} \label{LYZp} \emph{(Lutwak et.\ al \textbf{\cite{LYZ2006}})} If $1 \leq p < \infty$ and $f \in W^{1,p}(\mathbb{R}^n)$ is not $0$ a.e., then there exists a unique origin-symmetric convex body $\langle f \rangle_p$ with nonempty interior such that
\[\int_{\mathbb{R}^n} g(\nabla f(x))^p\,dx = \frac{1}{\left | \langle f \rangle_p \right |} \int_{\mathbb{S}^{n-1}} g(u)^p\,dS_p\left ( \langle f \rangle_p,u \right )  \]
for every even continuous function $g: \mathbb{R}^n \rightarrow [0,\infty)$ that is positively $1$-homogeneous.
\end{theorem}

For $p = 1$, where no volume-normalization is necessary, Wang \textbf{\cite{wang12}} extended Theorem \ref{LYZp} to functions of bounded variation. His result, which we state in the following, is the critical ingredient to deduce our Theorem \ref{mainthm3} from Theorem \ref{mainthm1}.

\begin{theorem} \label{LYZbv} \emph{(Wang \textbf{\cite{wang12}})} If $f \in BV(\mathbb{R}^n)$ is not $0$ a.e., then there exists a unique origin-symmetric convex body $\langle f \rangle$ with nonempty interior such that
\[\int_{\mathbb{R}^n} g(\sigma_f(x))\,d|Df|(x) = \int_{\mathbb{S}^{n-1}} g(u)\,dS\left ( \langle f \rangle,u \right )  \]
for every even continuous function $g: \mathbb{R}^n \rightarrow \mathbb{R}$ that is positively $1$-homogeneous.
\end{theorem}

The convex bodies $\langle f \rangle_p$ encode the geometry of the level sets of $f$. More precisely, if $f$ coincides up to translation a.e.\ on $\mathbb{R}^n$ with the convex symmetrization $f^K$ with respect to
an origin-symmetric $K \in \mathcal{K}^n$, then $\langle f \rangle_p$ is a dilate of $K$. In particular, the following proposition holds.

\begin{prop} \label{propertLYZ} \emph{(\!\!\textbf{\cite{LYZ2006, wang12}})} Let $K \in \mathcal{K}^n$ be an origin-symmetric convex body with nonempty interior.
\begin{enumerate}
\item[(i)] $\langle 1_K \rangle = K$.
\item[(ii)] For every $1 < p < \infty$, there exists $f \in W^{1,p}(\mathbb{R}^n)$ such that $\langle f \rangle_p = K$.
\end{enumerate}
\end{prop}

\pagebreak

\vspace{1cm}

\centerline{\large{\bf{ \setcounter{abschnitt}{4}
\arabic{abschnitt}. Auxiliary results}}}

\reseteqn \alpheqn \setcounter{theorem}{0}

\vspace{0.6cm}

In the following we prove an important relation between the $L_p$ projection body operator $\Pi_p$ and a family of $L_p$ Minkowski valuations which extend the operators $\Phi^{\mu}$. At the end of the section we recall Sobolev type inequalities for the operators $\langle f \rangle$ and $\langle f \rangle_p$ due to Wang \textbf{\cite{wang12}} and Lutwak, Yang, and Zhang \textbf{\cite{LYZ2006}}, respectively.

\vspace{0.2cm}

We begin by briefly recalling the one-to-one correspondence between measures on the homogeneous space $\mathbb{S}^{n-1} \cong \mathrm{SO}(n)/\mathrm{SO}(n-1)$ and right $\mathrm{SO}(n-1)$ invariant measures on $\mathrm{SO}(n)$.
To this end suppose that $\mu$ is a measure on $\mathbb{S}^{n-1}$. Since for every $u \in \mathbb{S}^{n-1}$, there exists $\vartheta_u \in \mathrm{SO}(n)$ such that $\vartheta_u \bar{e} = u$, it follows from Riesz's representation theorem that there exists a (unique) measure $\breve{\mu}$ on $\mathrm{SO}(n)$ such that
\[\int_{\mathrm{SO}(n)} g(\phi)\,d\breve{\mu}(\phi)= \int_{\mathbb{S}^{n-1}} \int_{\mathrm{SO}(n-1)}g(\vartheta_u \theta)\,d\theta\,d\mu(u) \]
for every $g \in C(\mathrm{SO}(n))$. Here, the inner integration on the right hand side is with respect to the Haar probability measure on $\mathrm{SO}(n-1)$ (which implies that the right hand side is independent of the choice of $\vartheta_u$). Moreover, it is easy to see that $\breve{\mu}$ is right $\mathrm{SO}(n-1)$ invariant and that
\begin{equation} \label{mulift}
\int_{\mathbb{S}^{n-1}} f(w)\,d\mu(w) = \int_{\mathrm{SO}(n)} f(\phi \bar{e})\,d\breve{\mu}(\phi)
\end{equation}
for every $f \in C(\mathbb{S}^{n-1})$. In other words, the pushforward of $\breve{\mu}$ under the natural projection $\pi: \mathrm{SO}(n) \rightarrow \mathbb{S}^{n-1}$, $\pi(\vartheta)=\vartheta \bar{e}$, is given by $\mu$. For more details, see \textbf{\cite{Schu09}}.

\vspace{0.2cm}

For the proof of Theorem \ref{mainthm4}, we require an $L_p$ extension of Theorem \ref{mainthm1}. To this end, we first have to define $L_p$ analogues of the operators $\Phi^{\mu}$.

\vspace{0.3cm}

\noindent {\bf Definition} \emph{Suppose that $1 \leq p < \infty$ and that $\mu$ is an even, zonal measure on $\mathbb{S}^{n-1}$. Let $Z_p^{\mu}(\bar{e})$ be the zonoid generated by $\mu$.
For a convex body $K \in \mathcal{K}^n$ containing the origin in its interior, we define the convex body $\Phi^{\mu}_pK \in \mathcal{K}^n$ by}
\begin{equation} \label{phiplpsurf}
h(\Phi_p^{\mu} K,u)^p = \int_{\mathbb{S}^{n-1}} h(Z_p^{\mu}(u),v)^p\,dS_p(K,v), \qquad u \in \mathbb{S}^{n-1}.
\end{equation}

\vspace{0.3cm}

Comparing definitions (\ref{phiplpsurf}) and (\ref{piplpsurfmeas}), we see that the operators $\Phi_p^{\mu}$ and $\Pi_p$ are related in the same way as $\Phi^{\mu}$ and $\Pi$. In particular, when $\mu$ is discrete, $\Phi^{\mu}_p$ is a multiple of $\Pi_p$.

With our next lemma, we will make the relation between $\Phi_p^{\mu}$ and $\Pi_p$ even more precise. To this end assume that $\mu$ is a \emph{zonal} measure on $\mathbb{S}^{n-1}$ and define a measure $\breve{\mu}_u$ on $\mathrm{SO}(n)$ as the pushforward of $\breve{\mu}$ under the conjugation map $c_u(\phi) = \vartheta_u \phi \vartheta_u^{-1}$.
Note that the $\mathrm{SO}(n-1)$ invariance of $\mu$ implies that $\breve{\mu}$ is $\mathrm{SO}(n-1)$ bi-invariant and, thus, the measure $\breve{\mu}_u$ is well-defined (that is, independent of the choice of $\vartheta_u$).

\begin{lem} \label{critrepLp} If $\mu$ is an even, zonal measure on $\mathbb{S}^{n-1}$ and $1  \leq p < \infty$, then
\[h(\Phi_p^{\mu}K,u)^p = \frac{1}{a_{n,p}} \int_{\mathrm{SO}(n)} h(\Pi_p K,\phi u)^p\,d\breve{\mu}_u(\phi), \qquad u \in \mathbb{S}^{n-1},  \]
for every $K \in \mathcal{K}^n$ containing the origin in its interior.
\end{lem}
{\it Proof.} For the $L_p$ zonoid generated by $\mu$ we have, by (\ref{zonoidmeas}) and (\ref{mulift}),
\[h(Z^{\mu}_p(\bar{e}),v)^p = \int_{\mathbb{S}^{n-1}} |v \cdot w|^p\,d\mu(w) = \int_{\mathrm{SO}(n)} |v \cdot \phi\bar{e}|^p\,d\breve{\mu}(\phi)  \]
for every $v \in \mathbb{S}^{n-1}$. Letting $u = \vartheta_{u}\bar{e}$ and using (\ref{sptfctcomprot}), we therefore obtain
\[h(Z^{\mu}_p(u),v)^p = h(Z^{\mu}_p(\bar{e}),\vartheta_u^{-1}v)^p = \int_{\mathrm{SO}(n)}\!\!\!\! |v \cdot \vartheta_u\phi \vartheta_u^{-1}u|^p\,d\breve{\mu}(\phi)
=  \int_{\mathrm{SO}(n)}\!\!\!\! |v \cdot \phi u|^p\,d\breve{\mu}_u(\phi)  \]
for every $u, v \in \mathbb{S}^{n-1}$. Thus, by (\ref{phiplpsurf}), Fubini's theorem, and (\ref{piplpsurfmeas}), we arrive at the desired relation
\[h(\Phi_p^{\mu}K,u)^p = \int_{\mathbb{S}^{n-1}}  \int_{\mathrm{SO}(n)}\!\!\!\! |v \cdot \phi u|^p\,d\breve{\mu}_u(\phi)\,dS_p(K,v) = \frac{1}{a_{n,p}} \int_{\mathrm{SO}(n)}\!\!\!\! h(\Pi_p K,\phi u)^p\,d\breve{\mu}_u(\phi).  \]
\hfill $\blacksquare$

\vspace{0.3cm}

Combining (\ref{pi1pi}) with the special case $p = 1$ of Lemma \ref{critrepLp}, we obtain the following:

\begin{koro} \label{critrep} If $\mu$ is an even, zonal measure on $\mathbb{S}^{n-1}$, then
\[h(\Phi^{\mu}K,u) = 2 \int_{\mathrm{SO}(n)} h(\Pi K,\phi u)\,d\breve{\mu}_u(\phi), \qquad u \in \mathbb{S}^{n-1},\]
for every $K \in \mathcal{K}^n$.
\end{koro}

Corollary \ref{critrep} is a key ingredient in the proof of Theorem \ref{mainthm2}. We also remark that if $\Phi: \mathcal{K}^n \rightarrow \mathcal{K}^n$ is any even, continuous Minkowski valuation of degree $n - 1$ which is $\mathrm{SO}(n)$ equivariant and translation invariant, then a similar relation to the projection body exists, however not necessarily with a non-negative measure $\breve{\mu}_u$ but a distribution. This follows from a characterization of such Minkowski valuations by the second author (see \textbf{\cite{Schu06a}} and the discussion at the end of Section 6).

\vspace{0.2cm}

As we shall see in Section 6 (and was previously noted in \textbf{\cite{LYZ2006}} and \textbf{\cite{wang12}}), the convexification operators $\langle f \rangle_p$ and $\langle g \rangle$ for $f \in W^{1,p}(\mathbb{R}^n)$ and $g \in BV(\mathbb{R}^n)$, respectively, allow for effortless proofs of Sobolev type inequalities once the corresponding isoperimetric inequalities have been established. This is due to the following sharp volume inequalities.

\begin{theorem} \label{volineqLYZop} \emph{(Wang \textbf{\cite{wang12}})}  If $f \in BV(\mathbb{R}^n)$, then
\[|\langle f \rangle |^{(n-1)/n} \geq \|f\|_{\frac{n}{n-1}}  \]
with equality if and only if $f$ is a multiple of the characteristic function of a centrally symmetric convex body.
\end{theorem}

The $L_p$ analogue of Theorem \ref{volineqLYZop} for Sobolev functions is the content of the following earlier result of Lutwak, Yang, and Zhang \textbf{\cite{LYZ2006}}.

\begin{theorem} \label{volineqLYZpop} \emph{(Lutwak et al.\ \textbf{\cite{LYZ2006}})}  If $1 < p < n$ and $f \in W^{1,p}(\mathbb{R}^n)$, then
\begin{equation} \label{lyzpabsch}
|\langle f \rangle_p |^{-1/n} \geq \tilde{c}_{n,p} \|f\|_{p^*},
\end{equation}
where the optimal constant $\tilde{c}_{n,p}$ is given by
\[ \tilde{c}_{n,p} =\left (\frac{n-p}{p-1}\right )^{1-1/p} \left ( \frac{\Gamma(\textstyle{\frac{n}{p}})\Gamma(n+1-\textstyle{\frac{n}{p}})}{\Gamma(n)}\right )^{1/n}.  \]
\end{theorem}

An infinite-dimensional set of extremal functions for inequality (\ref{lyzpabsch}) was exhibited in \textbf{\cite{LYZ2006}}. However, as far as we know, the equality conditions of (\ref{lyzpabsch}) have not been completely settled yet.

\vspace{1cm}

\centerline{\large{\bf{ \setcounter{abschnitt}{5}
\arabic{abschnitt}. Proof of the main results}}}

\reseteqn \alpheqn \setcounter{theorem}{0}

\vspace{0.6cm}

We are now in a position to prove Theorems \ref{mainthm1} and \ref{mainthm2} as well as their $L_p$ versions.
Their proofs are based largely on the following $L_p$ extension of the right hand inequality of Theorem \ref{mainthm2} (which, in turn, uses Lemma \ref{critrepLp}).

\begin{theorem} \label{maininequ17Lp} Suppose that $1 \leq p < \infty$. If $\mu$ is an even, zonal measure on $\mathbb{S}^{n-1}$ such that $\mu(\mathbb{S}^{n-1}) = a_{n,p}$ and
$K \in \mathcal{K}^n$ contains the origin in its interior, then
\begin{equation} \label{mainineqLp}
|\Phi^{\mu,\circ}_pK| \leq |\Pi^\circ_p K|
\end{equation}
with equality if and only if $\mu$ is discrete or $\Pi_p K$ is a Euclidean ball.
\end{theorem}

\noindent {\it Proof.}  First note that the normalization $\mu(\mathbb{S}^{n-1})=a_{n,p}$, with $a_{n,p}$ given by (\ref{piplpsurfmeas}), was chosen such that there is equality in (\ref{mainineqLp}) when $\mu$ is discrete.

In order to establish (\ref{mainineqLp}), we use the polar coordinate formula for volume (\ref{polarcoordvol}) and Lemma \ref{critrepLp} to obtain
\[|\Phi^{\mu,\circ}_pK| = \frac{1}{n} \int_{\mathbb{S}^{n-1}}\!\!\!\! h(\Phi^{\mu}_pK,u)^{-n}\,du =
\frac{1}{n} \int_{\mathbb{S}^{n-1}}\! \left ( \frac{1}{a_{n,p}} \int_{\mathrm{SO}(n)}\!\!\!\!\! h(\Pi_p K,\phi u)^p\,d\breve{\mu}_u(\phi)  \right )^{-n/p}du.  \]
Since $\mu(\mathbb{S}^{n-1}) = a_{n,p}$ implies $\breve{\mu}_u(\mathrm{SO}(n)) = a_{n,p}$, Jensen's inequality yields
\begin{equation} \label{jensencritLp}
|\Phi^{\mu,\circ}_pK| \leq \frac{1}{na_{n,p}} \int_{\mathbb{S}^{n-1}} \int_{\mathrm{SO}(n)} h(\Pi_p K,\phi u)^{-n}\,d\breve{\mu}_u(\phi)\,du.
\end{equation}
If follows from definitions (\ref{phiplpsurf}) and (\ref{piplpsurfmeas}) that $\Phi^{\mu}_p$ and $\Pi_p$ are $\mathrm{SO}(n)$ equivariant. Hence, replacing
$K$ by $\theta K$ in (\ref{jensencritLp}) and also using the $\mathrm{SO}(n)$ equivariance of the polar map and (\ref{sptfctcomprot}), we arrive at
\[ |\Phi^{\mu,\circ}_pK| \leq \frac{1}{na_{n,p}} \int_{\mathbb{S}^{n-1}} \int_{\mathrm{SO}(n)} \!\!\! h(\Pi_p K,\theta^{-1}\phi u)^{-n}\,d\breve{\mu}_u(\phi)\,du.  \]
Integrating now both sides with respect to the Haar probability measure on $\mathrm{SO}(n)$, followed by Fubini's theorem, we conclude that
\begin{equation} \label{proof17aLp}
 |\Phi^{\mu,\circ}_pK| \leq \frac{1}{na_{n,p}} \int_{\mathbb{S}^{n-1}} \int_{\mathrm{SO}(n)} \int_{\mathrm{SO}(n)}\!\!\! h(\Pi_p K,\theta^{-1}\phi u)^{-n}\,d\theta\,d\breve{\mu}_u(\phi)\,du.
\end{equation}
The invariance of the Haar measure and the fact that $\breve{\mu}_u(\mathrm{SO}(n)) = a_{n,p}$ imply that the right-hand side of
inequality (\ref{proof17aLp}) is equal to
\[ \frac{1}{n} \int_{\mathbb{S}^{n-1}} \int_{\mathrm{SO}(n)} \!\!\! h(\Pi_p K,\theta u)^{-n}\,d\theta\,du. \]
Using this in (\ref{proof17aLp}) followed by Fubini's theorem, we conclude from (\ref{sptfctcomprot}) that
\[|\Phi^{\mu,\circ}_pK| \leq \frac{1}{n} \int_{\mathrm{SO}(n)} \int_{\mathbb{S}^{n-1}} \!\!\! h(\Pi_p K,\theta u)^{-n}\,du\,d\theta
= \int_{\mathrm{SO}(n)} |\theta^{-1} \Pi^\circ_p K| \,d\theta = |\Pi^\circ_p K|. \]

By the above arguments, equality holds in (\ref{mainineqLp}) if and only if we have equality in (\ref{jensencritLp}). By the equality condition of Jensen's inequality this is the case if and only if
for every $u \in \mathbb{S}^{n-1}$ there exists $c_{u} \in \mathbb{R}^+$ such that
\[h(\Pi_p K,\phi u)=c_{u} \mbox{ for } \breve{\mu}_u\mbox{-a.e.\ } \phi \in \mathrm{SO}(n).  \]
But since $\breve{\mu}_u$ is the pushforward of $\breve{\mu}$ under the conjugation $c_u$ and $\mu$ is the pushforward of $\breve{\mu}$ under the projection $\pi$ (see Section 4), we see that equality holds in (\ref{jensencritLp}) if and only if for every $u \in \mathbb{S}^{n-1}$ and every $\vartheta \in \mathrm{SO}(n)$ such that $\vartheta\bar{e}=u$ there exists $c_{\vartheta_u} \in \mathbb{R}^+$ such that $h(\Pi_p K,\vartheta_uv) = c_{\vartheta_u}$ for $\mu$-a.e.\ $v \in \mathbb{S}^{n-1}$. Clearly, this is the case if and only if for every $\vartheta \in \mathrm{SO}(n)$, there exist $c_{\vartheta} \in \mathbb{R}^+$ such that
\begin{equation} \label{equality1Lp}
h(\Pi_p K,\vartheta v)=c_{\vartheta} \mbox{ for } \mu\mbox{-a.e.\ } v \in \mathbb{S}^{n-1}.
\end{equation}
We want to show that if $\mu$ is \emph{not discrete}, then (\ref{equality1Lp}) holds if and only if $h(\Pi_p K,\cdot)$ is constant on $\mathbb{S}^{n-1}$, or equivalently, if $\Pi_p K$ is a Euclidean ball. To this end, let $H_{w,t}=\{x \in \mathbb{R}^n: w \cdot x = t\}$ with $w \in \mathbb{S}^{n-1}$ and $t \in \mathbb{R}$ and note that since $\mu$ is zonal and not discrete, there exists $t_0 \in (-1,1)$ such that for all $\varepsilon > 0$ and all
$w \in H_{\bar{e},t_0} \cap \mathbb{S}^{n-1}$,
\[\mu(B_{\varepsilon}(w) \cap \mathbb{S}^{n-1}) > 0.  \]

\pagebreak

Now, by (\ref{equality1Lp}), we know that for every $\vartheta \in \mathrm{SO}(n)$ there is a set of $\mu$-measure zero $A_{\vartheta} \subseteq \mathbb{S}^{n-1}$ such that $h(\Pi_pK,\vartheta v) = c_{\vartheta}$ for all $v$ in the complement $A_{\vartheta}^{\mathrm{C}}$ of $A_{\vartheta}$. But for arbitrary $v \in H_{\bar{e},t_0} \cap \mathbb{S}^{n-1}$ and every $\varepsilon > 0$, we have that $B_{\varepsilon}(v) \cap A_{\vartheta}^{\mathrm{C}}$ is nonempty. Consequently, there is a sequence $v_k \in A_{\vartheta}^{\mathrm{C}}$ converging to $v$, which, by (\ref{equality1Lp}) and the continuity of support functions, shows that for every $\vartheta \in \mathrm{SO}(n)$ there exist $c_{\vartheta} \in \mathbb{R}^+$ such that
\begin{equation} \label{equalitychristoph}
h(\Pi_p K,\vartheta v)=c_{\vartheta} \mbox{ for all } v \in H_{\bar{e},t_0} \cap \mathbb{S}^{n-1}.
\end{equation}
In particular,
\[h(\Pi_p K,v) = c_{\mathrm{id}} \mbox{ for all } v \in H_{\bar{e},t_0} \cap \mathbb{S}^{n-1}.\]
Hence, applying a rotation $\theta \in \mathrm{SO}(n)$ to this equation and using (\ref{equalitychristoph}) again, yields
\[h(\Pi_p K,u) = c_{\theta} \mbox{ for all } u \in H_{\theta\bar{e},t_0} \cap \mathbb{S}^{n-1}.  \]
Now, if we choose $\theta$ such that $H_{\bar{e},t_0} \cap \mathbb{S}^{n-1}$ and $H_{\theta\bar{e},t_0} \cap \mathbb{S}^{n-1}$ intersect, then it follows that $c_{\mathrm{id}} = c_{\theta}$. Since we can reach any point on $\mathbb{S}^{n-1}$ by finitely many iterations of this argument, we obtain that $h(\Pi_p K,v) = c_{\mathrm{id}}$ for all $v \in \mathbb{S}^{n-1}$ as desired. \hfill $\blacksquare$

\vspace{0.4cm}

With the help of Theorem \ref{maininequ17Lp}, we can now complete the proof of Theorem \ref{mainthm2}.

\vspace{0.4cm}

\noindent {\it Proof of Theorem \ref{mainthm2}.} Let $K$ be a convex body in $\mathbb{R}^n$ with nonempty interior. First we want to show that
\begin{equation} \label{mainineqconvex}
\frac{n^n\omega_n^{n+1}}{\omega_{n-1}^n}|\partial K|^{-n} \leq |\Phi^{\mu,\circ}K| \leq |\Pi^\circ K|.
\end{equation}
Note that the normalization $\mu(\mathbb{S}^{n-1})=\frac{1}{2}$ ensures that there is equality in both inequalities of (\ref{mainineqconvex}) if $K$ is a Euclidean ball.
Moreoever, using (\ref{pi1pi}) and the translation invariance of both $\Phi^{\mu}$ and $\Pi$, we see that the right hand inequality of (\ref{mainineqconvex}) follows from
Theorem \ref{maininequ17Lp} by renormalization and that $|\Phi^{\mu,\circ}K| = |\Pi^\circ K|$ holds if and only if $\mu$ is discrete or $\Pi K$ is a Euclidean ball.

In order to establish the left hand inequality of (\ref{mainineqconvex}), we use the polar coordinate formula for volume (\ref{polarcoordvol}) and Jensen's inequality to obtain
\begin{equation} \label{leftineq1}
\left ( \frac{|\Phi^{\mu,\circ}K|}{\omega_n}  \right )^{-1/n} = \left (\frac{1}{n\omega_n}\int_{\mathbb{S}^{n-1}}\!\!\!\!\! h(\Phi^{\mu}K,u)^{-n}\,du \right )^{-1/n} \leq \frac{1}{n\omega_n}\int_{\mathbb{S}^{n-1}}
 \!\!\!\!\! h(\Phi^{\mu}K,u)\,du.
\end{equation}
From (\ref{phisurfmeas}) and the remark following it, Fubini's theorem, and (\ref{meanwidthzonoid}), we conclude
\[\int_{\mathbb{S}^{n-1}} \!\!\! h(\Phi^{\mu}K,u)\,du = \int_{\mathbb{S}^{n-1}}\int_{\mathbb{S}^{n-1}} h(Z^{\mu}(v),u)\,du\,dS(K,v) = \omega_{n-1}|\partial K|  \]
which yields the desired inequality. Note that equality holds in (\ref{leftineq1}), and thus in the left hand inequality of (\ref{mainineqconvex}), if and only if $h(\Phi^{\mu}K,\cdot)$ is constant, that is, if $\Phi^{\mu}K$ is a ball.

In order to extend the statement to a body $L \subseteq \mathbb{R}^n$ of finite perimeter with non-empty interior, note that by (\ref{phisurfmeas}), Theorem \ref{LYZbv}, and (\ref{Dvonchar}), we have for $u \in \mathbb{S}^{n-1}$,
\[h(\Phi^{\mu}\langle 1_L \rangle,u) = \int_{\mathbb{S}^{n-1}}\!\!\! h(Z^{\mu}(u),v)\,dS(\langle 1_L \rangle,v) = \int_{\partial^*L}\!\! h(Z^{\mu}(u),\nu_L(x))\,dx = h(\Phi^{\mu}L,u),   \]
that is, $\Phi^{\mu}L = \Phi^{\mu}\langle 1_L \rangle$. In particular, we also have $|\partial^*L| = |\partial \langle 1_L \rangle|$ and $\Pi L = \Pi\langle 1_L \rangle$. Thus, the desired inequalities along with their equality conditions follow from the first part of the proof.
\hfill $\blacksquare$

\vspace{0.3cm}

Theorem \ref{mainthm1} can now be deduced easily by combining the right-hand inequality of Theorem \ref{mainthm2}  with the Petty projection inequality.

\vspace{0.3cm}

\noindent {\it Proof of Theorem \ref{mainthm1}.} By our convention that $\infty \cdot 0 = 0$, we may assume that $L$ has nonempty interior. Moreover, we may also assume (by rescaling if necessary) that $\mu(\mathbb{S}^{n-1}) = \frac{1}{2}$ and that $\mu$ is not discrete. Hence, combining Theorem \ref{mainthm2} with the Petty projection inequality yields
\[|\Phi^{\mu,\circ}L||L|^{n-1} \leq |\Pi^{\circ}L||L|^{n-1} \leq \frac{\omega_n^n}{\omega_{n-1}^n},  \]
where equality holds in the left-hand inequality if and only if $\Pi L$ is a Euclidean ball and equality holds in the right-hand inequality if and only if $L$ is an ellipsoid up to a Lebesgue null set.

Now since
\[\Pi(AL) = |\mathrm{det}\, A|\,A^{-\mathrm{T}} \Pi L  \]
for any $A \in \mathrm{GL}(n)$ (cf.\ \textbf{\cite[\textnormal{Theorem 4.1.5}]{gardner2ed}}), it follows that if $E$ is an ellipsoid, then $\Pi E$ is a ball if and only if $E$ is a Euclidean ball. Consequently, we obtain
\[ |\Phi^{\mu,\circ}L||L|^{n-1} \leq \frac{\omega_n^n}{\omega_{n-1}^n} \]
with equality if and only if $L$ is a Euclidean ball up to a Lebesgue null set. \hfill $\blacksquare$

\vspace{0.3cm}

By combining now Theorem \ref{maininequ17Lp} with the $L_p$ Petty projection inequality, we obtain an $L_p$ extension of Theorem \ref{mainthm1}.

\begin{theorem} \label{mainthm1Lp} Suppose that $1 < p < \infty$ and that $\mu$ is an even, zonal measure on $\mathbb{S}^{n-1}$. Among convex bodies $K \in \mathcal{K}^n$ containing the origin in its interior the volume product
\[|\Phi^{\mu,\circ}_pK|^p|K|^{n-p}   \]
is maximized by origin-symmetric Euclidean balls. If $\mu$ is not discrete, then origin-symmetric Euclidean balls are the only maximizers. If $\mu$ is discrete, then
$K$ is a maximizer if and only if it is an origin-symmetric ellipsoid.
\end{theorem}

\noindent {\it Proof.} Let $K \in \mathcal{K}^n$ contain the origin in its interior and assume without loss of generality that $\mu(\mathbb{S}^{n-1}) = a_{n,p}$, where $a_{n,p}$ is given by (\ref{piplpsurfmeas}), and that $\mu$ is not discrete. Then, Theorem \ref{maininequ17Lp} and the $L_p$ Petty projection inequality yield
\[|\Phi^{\mu,\circ}_pK|^p|K|^{n-p} \leq |\Pi^{\circ}_pK|^p|K|^{n-p} \leq \omega_n^n  \]
with equality in the left-hand inequality if and only if $\Pi_p K$ is a Euclidean ball and equality in the right-hand inequality if and only if $K$ is an origin-centered ellipsoid. Since
\[\Pi_p(AK) = |\mathrm{det}\, A|^{1/p}\,A^{-\mathrm{T}} \Pi_p K  \]
for any $A \in \mathrm{GL}(n)$ (cf.\ \textbf{\cite{LYZ2000a}}), it follows that if $E$ is an origin-centered ellipsoid, then $\Pi_p E$ is a ball if and only if $E$ is a Euclidean ball. Consequently, we obtain
\[ |\Phi^{\mu,\circ}_pK|^p|K|^{n-p} \leq \omega_n^n \]
with equality if and only if $K$ is an origin-centered Euclidean ball. \hfill $\blacksquare$

\vspace{1cm}

\centerline{\large{\bf{ \setcounter{abschnitt}{6}
\arabic{abschnitt}. Applications to Sobolev inequalities}}}

\reseteqn \alpheqn \setcounter{theorem}{0}

\vspace{0.6cm}

We can now combine our isoperimetric inequalities from Theorems~\ref{mainthm1Lp}~and~\ref{mainthm1} with the volume inequalities for the operators $\langle f \rangle_p$ and $\langle f \rangle$, respectively, to deduce Theorems~\ref{mainthm4}~and~\ref{mainthm3}. Finally, we also relate our new Sobolev inequality from Theorem~\ref{mainthm3} to the sharp Gromov-Sobolev inequality
for general norms and conclude the paper with an open problem.

We begin with the proof of Theorem \ref{mainthm4} and, to this end, first specify the optimal constant $c_{n,p}$ appearing there:
\[c_{n,p} =\left (\frac{n-p}{p-1}\right )^{1-1/p} \left ( \frac{\Gamma(\textstyle{\frac{n}{p}})\Gamma(n+1-\textstyle{\frac{n}{p}})}{\Gamma(n+1)}\right )^{1/n}
\left (\frac{n\Gamma(\textstyle{\frac{n}{2}})\Gamma(\textstyle{\frac{p+1}{2}})}{\sqrt{\pi}\Gamma(\textstyle{\frac{n+p}{2}})}\right)^{1/p}.  \]

\vspace{0.3cm}

\noindent {\it Proof of Theorem \ref{mainthm4}.} We may assume that $f$ is not $0$ a.e.\ and, by rescaling if necessary, that $\mu(\mathbb{S}^{n-1}) = 1$.
Next, we recall that for every $K \in \mathcal{K}^n$ containing the origin in its interior and $\lambda > 0$, we have $S_p(\lambda K,\cdot) = \lambda^{n-p}S_p(K,\cdot)$.
Hence, taking
\[K = \left |\langle f \rangle_p\right |^{-\frac{1}{n-p}} \langle f \rangle_p,  \]
it follows from definition (\ref{phiplpsurf}) and Theorem \ref{LYZp} that for $u \in \mathbb{S}^{n-1}$,
\[h(\Phi_p^{\mu}K,u)^p = \frac{1}{\left |\langle f \rangle_p\right |} \int_{\mathbb{S}^{n-1}} h(Z_p^{\mu}(u),v)^p\,dS_p(\langle f \rangle_p,v) = \int_{\mathbb{R}^n} h(Z_p^{\mu}(u),\nabla f(x))^p\,dx.  \]
Since $h(Z_p^{\mu}(u),\nabla f(x))=\|\nabla f(x)\|_{Z^{\mu}_p(u)^\circ}$, this and the polar coordinate formula for volume (\ref{polarcoordvol}) imply that
\begin{equation} \label{prfThm4a}
|\Phi^{\mu,\circ}_pK|^{-1/n} = n^{1/n} \left ( \int_{\mathbb{S}^{n-1}} \left ( \int_{\mathbb{R}^n} \|\nabla f(x)\|^p_{Z^{\mu}_p(u)^\circ}\,dx \right )^{-n/p} du \right )^{-1/n}.
\end{equation}
Next, note that under the normalization $\mu(\mathbb{S}^{n-1})=1$ and our choice of the body $K$, Theorem \ref{mainthm1Lp} states that
\begin{equation} \label{prfThm4b}
|\Phi_p^{\mu,\circ}K|^{-1/n} \geq (a_{n,p}\omega_n)^{-1/p}|K|^{1/p^*} = (a_{n,p}\omega_n)^{-1/p} \left |\langle f \rangle_p \right |^{-1/n},
\end{equation}
where $a_{n,p}$ is the normalizing constant from (\ref{piplpsurfmeas}). Combining now (\ref{prfThm4a}), (\ref{prfThm4b}), and Theorem \ref{volineqLYZpop}, we obtain the desired inequality
\[\left ( \int_{\mathbb{S}^{n-1}} \left ( \int_{\mathbb{R}^n} \|\nabla f(x)\|^p_{Z^{\mu}_p(u)^\circ}\,dx \right )^{-n/p} du \right )^{-1/n} \geq \frac{n^{-1/n}}{(a_{n,p}\omega_n)^{1/p}} \left |\langle f \rangle_p \right |^{-1/n}
\geq c_{n,p}\|f\|_{p^*}.  \]

\vspace{-0.4cm}

\hfill $\blacksquare$

\vspace{0.3cm}

Before we turn to the proof of Theorem \ref{mainthm3}, we want to point out the following interesting special case of Theorem \ref{mainthm4} which is obtained by taking the measure $\mu$ to be
(uniformly) concentrated on $\mathbb{S}^{n-1} \cap \bar{e}^{\bot}$ with total mass $\mu(\mathbb{S}^{n-1})= a_{n-1,p}$.

\begin{koro} \label{sintranskoro} Suppose that $1 < p < n$. If $f \in W^{1,p}(\mathbb{R}^n)$, then
\[\left ( \int_{\mathbb{S}^{n-1}} \left ( \int_{\mathbb{R}^n} |\nabla f(x)|u^{\bot}|^p \,dx \right )^{-n/p} du \right )^{-1/n} \geq c_{n,p}a_{n-1,p}^{1/p} \|f\|_{p^*}.\]
\end{koro}

Kniefacz and the second author use in \textbf{\cite{KniefaczSchu1}} a refinement of the arguments from the proof of Theorem \ref{mainthm2} to generalize both the affine $L_p$ Sobolev inequality (\ref{affLpSob})
and Corollary \ref{sintranskoro} from one and $(n-1)$-dimensional projection averages, respectively, to averages over projections onto $k$-dimensional subspaces.

\vspace{0.3cm}

\noindent {\it Proof of Theorem \ref{mainthm3}.} We may assume that $f$ is not $0$ a.e., $\mu$ is not discrete and, by rescaling if necessary, that $\mu(\mathbb{S}^{n-1}) = 1$.
Then, by (\ref{phisurfmeas}) and Theorem \ref{LYZbv}, we have for $u \in \mathbb{S}^{n-1}$,
\[h(\Phi^{\mu}\langle f \rangle,u) = \int_{\mathbb{S}^{n-1}} h(Z^{\mu}(u),v)\,dS(\langle f \rangle,v) = \int_{\mathbb{R}^n} h(Z^{\mu}(u),\sigma_f(x))\,d|Df|(x).  \]
Since $h(Z^{\mu}(u),\sigma_f(x))=\|\sigma_f(x)\|_{Z^{\mu}(u)^\circ}$, this and the polar coordinate formula for volume (\ref{polarcoordvol}) imply that
\begin{equation} \label{prfThm3a}
|\Phi^{\mu,\circ}\langle f \rangle|^{-1/n} = n^{1/n} \left ( \int_{\mathbb{S}^{n-1}} \left ( \int_{\mathbb{R}^n} \|\sigma_f\|_{Z^{\mu}(u)^\circ}\,d|Df| \right )^{-n} du \right )^{-1/n}.
\end{equation}
Next, note that under the normalization $\mu(\mathbb{S}^{n-1})=1$, Theorem \ref{mainthm1} yields
\begin{equation} \label{prfThm3b}
|\Phi^{\mu,\circ}\langle f \rangle|^{-1/n} \geq \frac{2\omega_{n-1}}{\omega_n} \left |\langle f \rangle \right |^{(n-1)/n},
\end{equation}
with equality if and only if $\langle f \rangle$ is a Euclidean ball. Combining now (\ref{prfThm3a}), (\ref{prfThm3b}), and Theorem \ref{volineqLYZop}, we obtain the desired inequality
\begin{equation} \label{prfThm3c}
\left ( \int_{\mathbb{S}^{n-1}} \left ( \int_{\mathbb{R}^n} \|\sigma_f\|_{Z^{\mu}(u)^\circ}\,d|Df| \right )^{-n} du \right )^{-1/n} \geq \frac{2\omega_{n-1}}{n^{1/n}\omega_n}\, \|f\|_{\frac{n}{n-1}}.
\end{equation}
Moreover, from the equality conditions of Theorem \ref{volineqLYZop} and (\ref{prfThm3b}) together with Proposition \ref{propertLYZ} (i), we see that equality holds in (\ref{prfThm3c}) if and only if
$f$ is a multiple of the characteristic function of a Euclidean ball in $\mathbb{R}^n$.
\hfill $\blacksquare$

\vspace{0.3cm}

Choosing again the measure $\mu$ concentrated on $\mathbb{S}^{n-1} \cap \bar{e}^{\bot}$ with total mass $\mu(\mathbb{S}^{n-1})=2\omega_{n-2}$,
we obtain the following corollary as a special case of Theorem \ref{mainthm3}:

\begin{koro} If $f \in BV(\mathbb{R}^n)$, then
\[\left ( \int_{\mathbb{S}^{n-1}} \left ( \int_{\mathbb{R}^n} |\sigma_f(x)|u^{\bot}| \,d|Df|(x) \right )^{-n} du \right )^{-1/n} \geq  \frac{4\omega_{n-1}\omega_{n-2}}{n^{1/n}\omega_n}\, \|f\|_{\frac{n}{n-1}}\]
with equality if and only if $f$ is a multiple of the characteristic function of a Euclidean ball in $\mathbb{R}^n$.
\end{koro}

Let us also remark here that the above proofs of Theorems \ref{mainthm4} and \ref{mainthm3} combined with Theorems \ref{mainthm2} and \ref{mainthm1Lp} show that among the established families of Sobolev inequalities the only affine invariant ones, the affine $L_p$ Sobolev inequality and the affine Sobolev-Zhang inequality, respectively, are the strongest ones.

\vspace{0.2cm}

We conclude the paper by relating our new family of Sobolev inequalities from Theorem \ref{mainthm3} to the Gromov-Sobolev inequality for general norms (cf.\ \textbf{\cite[\textnormal{Appendix}]{milmschechtm}}). Recall that the latter states that for any origin-symmetric $K \in \mathcal{K}^n$ with nonempty interior and $f \in BV(\mathbb{R}^n)$, we have
\begin{equation} \label{gromovsob}
\int_{\mathbb{R}^n} \|\sigma_f\|_{K^\circ}\,d|Df| \geq n|K|^{1/n}\|f\|_{\frac{n}{n-1}}
\end{equation}
with equality if and only if $f$ is a multiple of the characteristic function of a homothetic copy of $K$.

Now, if $\mu$ is an even, zonal measure on $\mathbb{S}^{n-1}$ (not concentrated on a great subsphere) and $Z^{\mu}(\bar{e})$ denotes again the zonoid generated by $\mu$, we can take
$K = Z^{\mu}(u)$ and average both sides of (\ref{gromovsob}) over all directions $u \in \mathbb{S}^{n-1}$ to obtain on one hand the inequality
\begin{equation} \label{grombettera}
\left ( \int_{\mathbb{S}^{n-1}} \left ( \int_{\mathbb{R}^n} \|\sigma_f\|_{Z^{\mu}(u)^\circ}\,d|Df| \right )^{-n} du \right )^{-1/n} \geq \frac{n^{(n-1)/n}|Z^{\mu}(\bar{e})|^{1/n}}{\omega_n^{1/n}}\, \|f\|_{\frac{n}{n-1}}.
\end{equation}
On the other hand, we can use (\ref{meanwidthzonoid}) to rewrite our Sobolev inequality from Theorem~\ref{mainthm3} as follows
\begin{equation} \label{grombetterb}
\left ( \int_{\mathbb{S}^{n-1}} \left ( \int_{\mathbb{R}^n} \|\sigma_f\|_{Z^{\mu}(u)^\circ}\,d|Df| \right )^{-n} du \right )^{-1/n} \geq \frac{n^{(n-1)/n}w(Z^{\mu}(\bar{e}))}{2}\, \|f\|_{\frac{n}{n-1}}.
\end{equation}
Comparing the right hand sides of (\ref{grombettera}) and (\ref{grombetterb}) using Urysohn's inequality (\ref{urysohn}),
it follows that our Sobolev inequality (\ref{grombetterb}) is stronger than the (simple) consequence (\ref{grombettera}) of the Gromov-Sobolev inequality.
Similar, considerations apply in the $L_p$~setting, when Theorem \ref{mainthm4} is compared with an average of the $L_p$~extension of the Gromov-Sobolev inequality by Cordero, Nazaret, and Villani \textbf{\cite{cordnazvill}}.

\vspace{0.2cm}

Note, however, that while inequality (\ref{grombettera}) remains true, when the zonoids $Z^{\mu}(u)$ are replaced by arbitrary convex bodies of revolution, it remains an open problem whether this is also possible for our inequality (\ref{grombetterb}). If this was the case, the corresponding isoperimetric inequality would generalize Theorem \ref{mainthm1} to \emph{all} even, continuous Minkowski valuations of degree $n-1$ which are $\mathrm{SO}(n)$ equivariant and
translation invariant.

\vspace{0.5cm}

\noindent {{\bf Acknowledgments} The second author was
supported by the European Research Council (ERC), Project number: 306445, and the Austrian Science Fund (FWF), Project number:
Y603-N26.

\begin{small}

\[ \begin{array}{ll} \mbox{Christoph Haberl} & \mbox{Franz Schuster} \\
\mbox{Vienna University of Technology \phantom{wwwwWW}} & \mbox{Vienna University of Technology} \\ \mbox{christoph.haberl@tuwien.ac.at} & \mbox{franz.schuster@tuwien.ac.at}
\end{array}\]

\end{small}


\begin{thebibliography}{99}
\footnotesize{
\parskip-0.1cm{

\bibitem{abardber11}
J. Abardia and A. Bernig, \emph{Projection bodies in complex vector spaces}, Adv. Math. {\bf 227} (2011), 830--846.

\bibitem{Alesker99}
S. Alesker, \emph{Continuous rotation invariant valuations on convex sets}, Ann. of Math. (2) {\bf 149} (1999), 977--1005.

\bibitem{Alesker01}
S. Alesker, \emph{Description of translation invariant valuations on convex sets with solution of P. McMullen's conjecture}, Geom. Funct. Anal. {\bf 11} (2001), 244--272.

\bibitem{ABS2011}
S. Alesker, A. Bernig, and F.E. Schuster, \emph{Harmonic analysis of translation invariant valuations}, Geom. Funct. Anal. {\bf 21} (2011), 751--773.

\bibitem{aubin76}
T. Aubin, \emph{Probl\`emes isop\'erim\'etriques et espaces de Sobolev},  J. Differential Geom. {\bf 11} (1976), 573--598.

\bibitem{BPSW2014}
A. Berg, L. Parapatits, F.E.Schuster, and M. Weberndorfer, \emph{Log-Concavity Properties of Minkowski Valuations}, arXiv:1411.7891.

\bibitem{bernigfu10}
A. Bernig and J.H.G. Fu, \emph{Hermitian integral geometry}, Ann. of Math. (2) {\bf 173} (2011), 907--945.

\bibitem{boeroezky2013}
K.J. B\"or\"oczky, \emph{Stronger versions of the Orlicz-Petty projection inequality}, J. Differential Geom. {\bf 95} (2013), 215–247.

\bibitem{BLYZ2013}
K.J. B\"or\"oczky, E. Lutwak, D. Yang, and G. Zhang, \emph{The logarithmic Minkowski problem}, J. Amer. Math. Soc. {\bf 26} (2013), 831--852.

\bibitem{Campi:Gronchi02a}
S. Campi and P. Gronchi, \emph{The $L_p$-Busemann-Petty centroid inequality}, Adv. Math. \textbf{167} (2002) 128--141.

\bibitem{chiachietal}
A. Cianchi, E. Lutwak, D. Yang, and G. Zhang, \emph{Affine Moser-Trudinger and Morrey-Sobolev inequalities}, Calc. Var. Partial Differential Equations {\bf 36} (2009), 419--436.

\bibitem{cordnazvill}
D. Cordero-Erausquin, B. Nazaret, and C. Villani, \emph{A mass-transportation approach to sharp Sobolev and Gagliardo-Nirenberg inequalities}, Adv. Math. {\bf 182} (2004), 307--332.

\bibitem{evansgariepy}
L.C. Evans and R.F. Gariepy, \emph{Measure theory and fine properties of functions}, Studies in Advanced Mathematics, CRC Press, Boca Raton, FL, 1992.

\bibitem{fedflem}
H. Federer and W. Fleming, \emph{Normal and integral currents}, Ann. Math. {\bf 72} (1960), 458--520.

\bibitem{gardner2ed}
R.J. Gardner, \emph{Geometric tomography}, Second ed., Cambridge University Press, Cambridge, 2006.

\bibitem{GooWei92}
P. Goodey and W. Weil, \emph{The determination of convex bodies from the mean of random sections}, Math. Proc. Cambridge Philos. Soc. {\bf 112} (1992), 419--430.

\bibitem{GooWei12}
P. Goodey and W. Weil, \emph{A uniqueness result for mean section bodies}, Adv. Math. {\bf 229} (2012), 596--601.

\bibitem{goodeyweil4}
P. Goodey and W. Weil, \emph{Sums of sections, surface area measures and the general Minkowski problem}, J. Differential Geom. {\bf 97} (2014), 477--514.

\bibitem{grinbergzhang99}
E. Grinberg and G. Zhang, \emph{Convolutions, transforms, and convex bodies}, Proc. London Math. Soc. (3) {\bf 78} (1999), 77--115.

\bibitem{haberl11}
C. Haberl, \emph{Minkowski valuations intertwining the special linear group}, J. Eur. Math. Soc. {\bf 14} (2012), 1565--1597.

\bibitem{habparap14}
C. Haberl and L. Parapatits, \emph{The centro-affine Hadwiger theorem}, J. Amer. Math. Soc. {\bf 27} (2014), 685--705.

\bibitem{habschu09}
C. Haberl and F.E. Schuster, \emph{General $L_p$ affine isoperimetric inequalities}, J. Differential Geom. {\bf 83} (2009), 1--26.

\bibitem{habschu09a}
C. Haberl and F.E. Schuster, \emph{Asymmetric affine $L_p$ Sobolev inequalities}, J. Funct. Anal. {\bf 257} (2009), 641--658.

\bibitem{habschuxia}
C. Haberl, F.E. Schuster, and J. Xiao, \emph{An asymmetric affine P\'olya-Szeg\H{o} principle}, Math. Ann. {\bf 352} (2012), 517--542.

\bibitem{haddad2016}
J. Haddad, C.H. Jim\'enez, and M. Montenegro, \emph{Sharp affine Sobolev type inequalities via the $L_p$ Busemann-Petty centroid inequality}, J. Funct. Anal. {\bf 271} (2016), 454--473.

\bibitem{kiderlen05}
M. Kiderlen, \emph{Blaschke- and Minkowski-Endomorphisms of convex bodies}, Trans. Amer. Math. Soc. {\bf 358} (2006), 5539--5564.

\bibitem{KniefaczSchu}
P. Kniefacz and F.E. Schuster, \emph{Affine vs. Euclidean rearrangement inequalities}, in preparation.

\bibitem{KniefaczSchu1}
P. Kniefacz and F.E. Schuster, \emph{Affine vs. Euclidean Sobolev inequalities}, in preparation.

\bibitem{Lin}
Y. Lin, \emph{Affine Orlicz P\'olya-Szeg\H{o} principle for log-concave functions}, J. Funct. Anal. {\bf 273} (2017), 3295--3326.

\bibitem{ludwig02} M. Ludwig, \emph{Projection bodies and valuations}, Adv. Math. {\bf 172} (2002), 158--168.

\bibitem{Ludwig:Minkowski}
M. Ludwig, \emph{Minkowski valuations}, Trans. Amer. Math. Soc. \textbf{357} (2005), 4191--4213.

\bibitem{Ludwig10a}
M. Ludwig, \emph{Minkowski areas and valuations}, J. Differential Geom. {\bf 86} (2010), 133--161.

\bibitem{Ludwig12}
M. Ludwig, \emph{Valuations on Sobolev spaces}, Amer. J. Math. {\bf 134} (2012), 827--842.

\bibitem{centro}
M. Ludwig and M. Reitzner, \emph{A classification of $\mathrm{SL}(n)$ invariant valuations}, Ann. of Math. (2) {\bf 172} (2010), 1223--1271.

\bibitem{LudXiaZhang}
M. Ludwig, J. Xiao, and G. Zhang, \emph{Sharp convex Lorentz-Sobolev inequalities}, Math. Ann. {\bf 350} (2011), 169--197.

\bibitem{lutwak84}
E. Lutwak, \emph{A general isepiphanic inequality}, Proc. Amer. Math. Soc. {\bf 90} (1984), 415--421.

\bibitem{lutwak93a}
E. Lutwak, \emph{The Brunn-Minkowski-Firey theory. I. Mixed volumes
  and the Minkowski problem}, J. Differential Geom. \textbf{38} (1993),
  131--150.

\bibitem{lutwak93}
E. Lutwak, \emph{Inequalities for mixed projection bodies}, Trans. Amer. Math. Soc. {\bf 339} (1993), 901--916.

\bibitem{lutwak96}
E. Lutwak, \emph{The Brunn-Minkowski-Firey theory. II. Affine and
  geominimal surface areas}, Adv. Math. \textbf{118} (1996), 244--294.

\bibitem{LYZ2000a}
E. Lutwak, D. Yang, and G. Zhang, \emph{$L_p$ affine isoperimetric inequalities}, J. Differential Geom. {\bf 56} (2000), 111--132.

\bibitem{LYZ2002}
E. Lutwak, D. Yang, and G. Zhang, \emph{Sharp affine $L_p$ Sobolev inequalities}, J. Differential Geom. {\bf 62} (2002), 17--38.

\bibitem{LYZ2004}
E. Lutwak, D. Yang, and G. Zhang, \emph{On the $L_p$-Minkowski problem}, Trans. Amer. Math. Soc. {\bf 356} (2004), 4359--4370.

\bibitem{LYZ2006}
E. Lutwak, D. Yang, and G. Zhang, \emph{Optimal Sobolev norms and the $L_p$ Minkowski problem}, Int. Math. Res. Not. 2006, Art. ID 62987, 21 pp.

\bibitem{LYZ2010a}
E. Lutwak, D. Yang, and G. Zhang, \emph{Orlicz projection bodies}, Adv. Math. {\bf 223} (2010), 220--242.

\bibitem{mareschschu}
G. Maresch and F.E. Schuster, \emph{The sine transform of isotropic measures}, Int. Math. Res. Not. IMRN 2012, 717--739.

\bibitem{mazya}
V.G. Maz'ya, \emph{Classes of domains and imbedding theorems for function spaces}, Dokl. Akad. Nauk. SSSR {\bf 133} (1960), 527--530.

\bibitem{milmschechtm}
V.D. Milman and G. Schechtman, \emph{Asymptotic theory of finite-dimensional normed spaces}, Lecture Notes in Mathematics, vol. 1200, Springer-Verlag, Berlin, 1986, With an appendix by M. Gromov.

\bibitem{nguyen}
V.H. Nguyen, \emph{New approach to the affine P\'olya-Szeg\H{o} principle and the stability version of the affine Sobolev inequality}, Adv. Math. {\bf 302} (2016), 1080--1110.

\bibitem{parapschu}
L. Parapatits and F.E. Schuster, \emph{The Steiner formula for Minkowski valuations}, Adv. Math. {\bf 230} (2012), 978--994.

\bibitem{parapwann}
L. Parapatits and T. Wannerer, \emph{On the inverse Klain map}, Duke Math. J. {\bf 162} (2013), 1895--1922.

\bibitem{parap2014a}
L. Parapatits, \emph{$\mathrm{SL}(n)$-contravariant $L_p$-Minkowski valuations}, Trans. Amer. Math. Soc. {\bf 366} (2014), 1195--1211.

\bibitem{parap2014b}
L. Parapatits, \emph{$\mathrm{SL}(n)$-covariant $L_p$-Minkowski valuations}, J. Lond. Math. Soc. {\bf 89} (2014), 397--414.

\bibitem{petty67}
C.M. Petty, \emph{Isoperimetric problems}, Proc. Conf. Convexity and Combinatorial Geometry (Univ. Oklahoma, 1971), University of Oklahoma, 1972, 26--41.

\bibitem{schneider74}
R. Schneider, \emph{Equivariant endomorphisms of the space of convex bodies}, Trans. Amer. Math. Soc. {\bf 194} (1974), 53--78.

\bibitem{schneider93}
R. Schneider, \emph{Convex Bodies: The Brunn--Minkowski Theory}, Second ed., Encyclopedia of Mathematics and its Applications 151, Cambridge University Press, Cambridge, 2013.

\bibitem{Schu06a}
F.E. Schuster, \emph{Convolutions and multiplier transformations of convex bodies}, Trans. Amer. Math. Soc. {\bf 359} (2007), 5567--5591.

\bibitem{Schu09}
F.E. Schuster, \emph{Crofton Measures and Minkowski Valuations}, Duke Math. J. {\bf 154} (2010), 1--30.

\bibitem{SchuWan11}
F.E. Schuster and T. Wannerer, \emph{$\mathrm{GL}(n)$ contravariant Minkowski valuations}, Trans. Amer. Math. Soc. {\bf 364} (2012), 815--826.

\bibitem{SchuWan13}
F.E. Schuster and T. Wannerer, \emph{Even Minkowski valuations},  Amer. J. Math. {\bf 137} (2015), 1651--1683.

\bibitem{SchuWan16}
F.E. Schuster and T. Wannerer, \emph{Minkowski valuations and generalized valuations}, J. Eur. Math. Soc., in press.

\bibitem{talenti76}
G. Talenti, \emph{Best constant in Sobolev inequality}, Ann. Math. Pura Appl. {\bf 110} (1976), 353--372.

\bibitem{wang12}
T. Wang, \emph{The affine Sobolev-Zhang inequality on $BV(\mathbb{R}^n)$}, Adv. Math. {\bf 230} (2012), 2457--2473.

\bibitem{wang13}
T. Wang, \emph{The affine P\'olya-Szeg\H{o} principle: equality cases and stability}, J. Funct. Anal. {\bf 265} (2013), 1728--1748.

\bibitem{wannerer10}
T. Wannerer, \emph{$\mathrm{GL}(n)$ equivariant Minkowski valuations}, Indiana Univ. Math. J. {\bf 60} (2011), 1655--1672.

\bibitem{zhang99}
G. Zhang, \emph{The affine Sobolev inequality}, J. Differential Geom. {\bf 53} (1999), 183--202.

}}
\end{thebibliography}
\end{document}